\documentclass[11pt,reqno]{amsart}
\usepackage[T1]{fontenc}
\usepackage{lmodern}
\usepackage{amssymb,mathtools}
\mathtoolsset{showonlyrefs=true}
\usepackage{microtype}
\usepackage{geometry}
\usepackage{booktabs}
\usepackage{graphicx}
\usepackage{hyperref}
\hypersetup{colorlinks=true,linkcolor=black,citecolor=black,urlcolor=blue,pdftitle={Postcritically Finite Unicritical Polynomials over CM Fields}}
\allowdisplaybreaks[2]
\newtheorem{theorem}{Theorem}[section]
\newtheorem{proposition}[theorem]{Proposition}
\newtheorem{lemma}[theorem]{Lemma}
\newtheorem{corollary}[theorem]{Corollary}
\newtheorem{remark}[theorem]{Remark}
\newcommand{\Q}{\mathbb{Q}}
\newcommand{\Z}{\mathbb{Z}}
\newcommand{\C}{\mathbb{C}}
\newcommand{\R}{\mathbb{R}}
\newcommand{\OO}{\mathcal{O}}
\DeclareMathOperator{\Gal}{Gal}
\DeclareMathOperator{\Cl}{Cl}
\DeclareMathOperator{\Nm}{N}
\newcommand{\PCF}{\mathord{\mathrm{PCF}}}
\DeclareMathOperator{\av}{Av}
\newcommand{\ab}{\mathrm{ab}}
\newcommand{\tr}{\mathrm{tr}}
\newcommand{\mui}{\mu_\infty}
\title[PCF unicritical polynomials over CM fields]{Postcritically Finite Unicritical Polynomials over CM Fields}
\author{Kaoru Sano}
\address{NTT Institute for Fundamental Mathematics, NTT Communication Science Laboratories, NTT, Inc., 2-4, Hikaridai, Seika-cho, Soraku-gun, Kyoto 619-0237, Japan}
\email{kaoru.sano@ntt.com}
\subjclass[2020]{Primary 37P15, 11R20\quad Secondary 11R18, 37F10}
\keywords{Postcritically finite polynomials, unicritical dynamics, CM fields, complex multiplication, algebraic integers, roots of unity, maximal abelian extension}
\begin{document}
    \begin{abstract}
        A rational map of degree at least two is postcritically finite
        (PCF) if every critical point has a finite forward orbit. For the unicritical family
        $f_{d,c}(z)=z^d+c$ with $d\ge2$, we call $c$ a PCF parameter if
        $f_{d,c}$ is PCF. As discussed by Baker and DeMarco, these
        parameters are analogous to $j$-invariants of elliptic curves
        with complex multiplication.
        For every $d\ge2$, we determine which PCF parameters lie in a
        CM field. The same list gives all PCF parameters over the maximal
        abelian extension of $\Q$.
        This extends the classifications of Noytaptim and Petsche over
        the maximal totally real extension of $\Q$ and of Ferraguti and
        Pagano over quadratic number fields. The proof uses the totally
        real classification in degree two, orbit bounds, Kronecker's
        theorem, Goksel's unit theorem, Lam--Leung's short-sum
        classification, cyclotomic averages, and known irreducibility
        results.
    \end{abstract}
    \maketitle

    \section{Introduction}\label{sec:introduction}
    Fix an integer $d\ge2$ and consider the unicritical family
    $f_{d,c}(z)=z^d+c$ with $c\in\C$. The only finite critical point is
    $0$, whose orbit starts $0\mapsto c\mapsto c^d+c\mapsto\cdots$.
    We call $f_{d,c}$ \emph{postcritically finite} (PCF), and $c$ a
    \emph{PCF parameter}, if this orbit is finite, equivalently eventually
    periodic. For example, $d=2$ and $c=-1$ give $0\mapsto-1\mapsto0$.
    Since infinity is fixed, this agrees with the definition of a PCF
    rational map, for which all critical orbits are finite.

    PCF polynomials play an important role in complex dynamics.
    They admit finite combinatorial descriptions and provide distinguished
    parameters in families. Critical portraits record the angles of selected
    external rays and classify marked PCF polynomials. This classification
    was first established for polynomials whose finite critical points are all
    strictly preperiodic \cite[Theorem~II]{BielefeldFisherHubbard}
    and then in general \cite[Theorems~1.1 and~1.3]{PoirierCriticalPortraits}.
    A Hubbard tree joins the critical orbits and records the dynamics,
    local degrees, and angles. Assigning external angles to a suitably
    extended tree gives a critical portrait
    \cite[Theorem~1.1 and Proposition~11.15]{PoirierHubbardTrees}.
    For PCF rational maps, Thurston equivalence gives uniqueness up to
    M\"obius conjugacy when the orbifold is hyperbolic
    \cite[Theorem~1]{DouadyHubbard}. Without that hypothesis, the only
    exceptions to uniqueness are flexible Latt\`es maps
    \cite[Theorem~6.2]{McMullenFamilies}.

    Each bounded hyperbolic component in the full space of monic centered
    polynomials of fixed degree has a unique PCF center
    \cite[Corollary~5.2, p.~196]{MilnorHyperbolic}.
    In the moduli space of degree-$d$ polynomials with marked critical
    points, the bifurcation measure has support equal to the closure of
    the locus where all finite critical points land on repelling cycles
    \cite[Theorem~9 and its proof]{DujardinFavre}.
    In the moduli space of degree-$d$ rational maps, the support is the
    closure of the PCF locus with no periodic critical points, after removing
    flexible Latt\`es maps from that locus \cite[Main Theorem, p.~493]{BuffEpstein}.
    Both closures are taken in the usual complex topology.

    For the family $f_{d,c}$, a finite critical orbit also imposes arithmetic
    conditions. A PCF parameter is an algebraic integer, and all its
    conjugates are PCF by Lemma~\ref{lem:escape}
    (see also \cite[Lemma~4.1]{FerragutiPagano}). More generally, PCF conjugacy
    classes of polynomials of fixed degree form a countable, Zariski-dense
    set of algebraic points in moduli space \cite[Proposition~2.6]{BakerDeMarco}.
    We ask which PCF parameters lie in a given field.
    Fix an algebraic closure $\overline{\Q}\subset\C$ and write
    $\PCF_d=\{c\in\overline{\Q}\mid f_{d,c}\text{ is PCF}\}$.

    \subsection*{Prior work and main results}
    Noytaptim and Petsche \cite{NP} classified PCF parameters
    in the maximal totally real extension $\Q^{\tr}$ for every $d\ge2$,
    with the following parameter lists.
    \begin{equation}\label{eq:NP}
        \PCF_d\cap\Q^{\tr}=
        \begin{cases}
            \{-2,-1,0\},&d=2,\\
            \{-1,0\},&d\ge4\text{ even},\\
            \{0\},&d\ge3\text{ odd}.
        \end{cases}
    \end{equation}
    Ferraguti and Pagano \cite[Theorem~4.2]{FerragutiPagano} classify
    the PCF parameters in quadratic number fields for every $d\ge2$.

    We extend both classifications to parameters in \emph{CM fields},
    the totally imaginary quadratic extensions of totally real number
    fields. For totally real $F$, the field $F(i)$ is CM, and every
    quadratic number field is totally real or CM. The containing CM
    field need not be Galois, and $\Q(c)$ need not be CM.

    Write $\zeta_N=e^{2\pi i/N}$ for $N\ge1$. For $d\ge3$ put
    \begin{equation}\label{eq:Sd}
        \mathcal{S}_d=
        \{0\}\cup\{c\in\overline{\Q}\mid c^{d-1}=-1\}
        \cup
        \begin{cases}
            \{c\in\overline{\Q}\mid c^{d-1}\in\{\zeta_3,\zeta_3^2\}\},
            &d\equiv0,1\pmod6,\\
            \varnothing,&\text{otherwise}.
        \end{cases}
    \end{equation}
    Each element of $\mathcal{S}_d$ is either zero or a root of unity.

    \begin{theorem}[Main theorem]\label{thm:main}
        For every integer $d\ge2$, we have the classification
        \begin{equation}\label{eq:main}
            \PCF_d\cap\bigcup_{F\text{ a CM field}}F
            =\begin{cases}
                \{0,-1,-2,i,-i\},& d=2,\\
                \mathcal{S}_d,& d\ge3.
            \end{cases}
        \end{equation}
    \end{theorem}

    Write $f_{d,c}^{\circ r}$ for the $r$th iterate, with
    $f_{d,c}^{\circ0}$ the identity. We say that $f_{d,c}$, or its
    parameter $c$, has \emph{exact critical orbit type} $(m,n)$ if
    $f_{d,c}^{\circ(m+n)}(0)=f_{d,c}^{\circ m}(0)$ and the first $m+n$
    orbit points, starting at $0$, are distinct. Here $m\ge0$ is the
    preperiod and $n\ge1$ the eventual period.
    Proposition~\ref{prop:high} and Corollary~\ref{cor:types} give the
    types of the listed parameters. The parameter counts are
    \begin{equation}\label{eq:cardinality}
        \#\left(\PCF_d\cap\bigcup_{F\text{ a CM field}}F\right)=
        \begin{cases}
            5,&d=2,\\
            3d-2,&d\ge3,\ d\equiv0,1\pmod6,\\
            d,&d\ge3,\ d\not\equiv0,1\pmod6.
        \end{cases}
    \end{equation}

\begin{figure}[htbp]
        \centering
        \includegraphics[width=\textwidth]{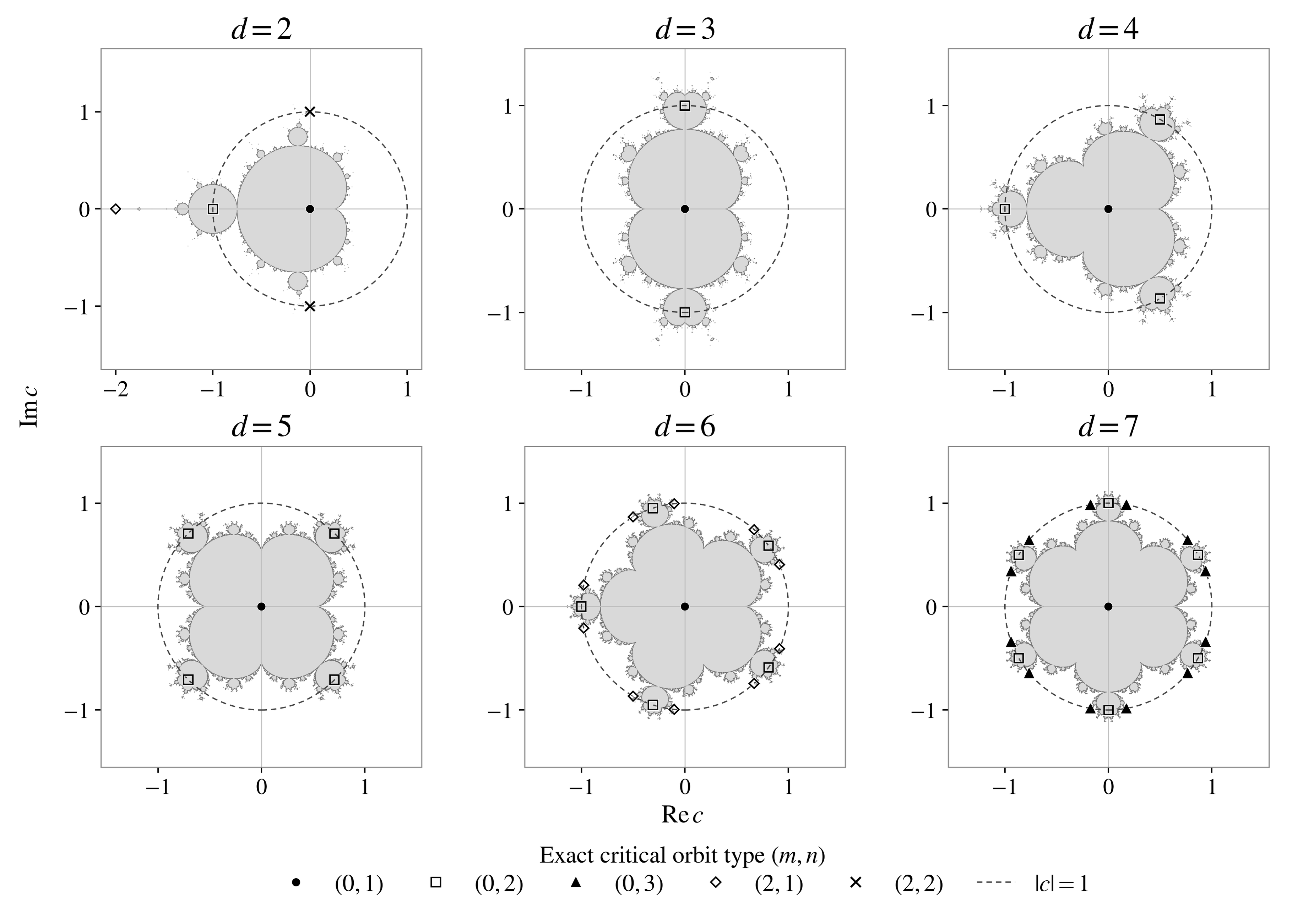}
        \caption{PCF parameters in CM fields for $d=2,\ldots,7$.
            Filled circles, open squares, filled triangles, open diamonds,
            and $\times$ marks indicate exact critical orbit types $(0,1)$, $(0,2)$,
            $(0,3)$, $(2,1)$, and $(2,2)$, respectively. The dashed circles indicate
            $|c|=1$. The shaded regions approximate the Multibrot sets
            $M_d=\{c\in\C\mid f_{d,c}^{\circ r}(0)\not\to\infty
            \text{ as }r\to\infty\}$.}
        \label{fig:cm-pcf-multibrot}
    \end{figure}

    Figure~\ref{fig:cm-pcf-multibrot} shows the classified parameters for
    $d=2,\ldots,7$ overlaid on the corresponding Multibrot sets, the sets
    of complex parameters for which the critical orbit of $0$ is bounded.

    Every parameter in \eqref{eq:main} is rational or a root of unity
    and hence lies in the maximal abelian extension $\Q^{\ab}$ of $\Q$.
    To see that the list contains every PCF parameter in $\Q^{\ab}$,
    we use the Kronecker--Weber theorem, for which Hilbert gave the
    first complete proof.
    \begin{theorem}[Kronecker--Weber {\cite{Kronecker1853,Weber1886,Weber1911}},
        Hilbert {\cite{Hilbert1896}}]\label{thm:KW}
        Every finite abelian extension of $\Q$ is contained in
        $\Q(\zeta_N)$ for some integer $N\ge1$. Consequently, we have
        \begin{equation}\label{eq:cyclotomic-union}
            \Q^{\ab}=\bigcup_{N\ge1}\Q(\zeta_N).
        \end{equation}
    \end{theorem}
    For $N\ge3$, $\Q(\zeta_N)$ is CM, and for $N=1,2$ it is contained
    in the CM field $\Q(\zeta_4)$. Thus every element of $\Q^{\ab}$ lies
    in a CM field, giving the following consequence of Theorem~\ref{thm:main}.

    \begin{corollary}[Maximal abelian extension]\label{cor:ab}
        For every $d\ge2$, the following two parameter sets coincide
        \begin{equation}\label{eq:abelian-cm-parameters}
            \PCF_d\cap\Q^{\ab}
            =\PCF_d\cap\bigcup_{F\text{ a CM field}}F.
        \end{equation}
        In particular, the lists in \eqref{eq:main} also classify PCF
        parameters over $\Q^{\ab}$.
    \end{corollary}

    In degree two, Ferraguti and Pagano's classification gives the five
    parameters in \eqref{eq:main} as the complete list over $\Q(i)$.
    Theorem~\ref{thm:main} and Corollary~\ref{cor:ab} show that this list
    remains complete in the union of all CM fields and over $\Q^{\ab}$.
    In degrees three and four, taking the union of the lists over all
    quadratic number fields likewise gives the complete list in each of
    these two settings. In every degree, the totally real classification
    gives the real parameters in $\PCF_d\cap\Q^{\ab}$, since every real
    element of $\Q^{\ab}$ is totally real.

    For comparison, let $K$ be a number field and $\mui$ the group
    of roots of unity in $\overline{\Q}$. Ingram and Mello
    \cite[Theorem~1]{IngramMello} prove finiteness of monic centered PCF
    polynomials of each degree $d\ge2$ over the cyclotomic closure $K(\mui)$.
    Zhang \cite[Theorems~1.1 and~1.2]{Zhang} proves a Bogomolov property
    for the critical height on polynomial moduli spaces over the maximal
    abelian extension $K^{\ab}$, yielding finitely many PCF conjugacy classes
    in each degree and hence finitely many unicritical PCF parameters,
    since each class contains at most $d-1$ maps of the form $f_{d,c}$.
    For $K=\Q$, both results give finiteness over $\Q^{\ab}$.
    Corollary~\ref{cor:ab} gives the explicit list in every degree for
    the unicritical family.

    The main step in degree two is to exclude strictly preperiodic
    critical orbits of eventual period at least three. Goksel's unit
    theorem \cite[Theorem~1.3(a)]{GokselOrbit} implies that $c$ and $c+1$
    are units for these orbit types. Their conjugates have absolute
    value at most $2$, so the CM unit decomposition stated below expresses
    each as a sum of two roots of unity. The two decompositions yield
    a vanishing sum of five roots of unity. Lam--Leung's theorem then
    restricts $c$ to four families (Lemma~\ref{lem:five}).
    We combine cyclotomic averages with direct orbit estimates to exclude
    these families for the stated orbit types (Lemma~\ref{lem:excludeHigher}).

    Benedetto and Goksel \cite{BenedettoGokselI} study the arithmetic of
    Misiurewicz parameters. Their sequel
    \cite[Proposition~2.1 and Corollary~1.8]{BenedettoGokselII} gives
    criteria for unit values of Misiurewicz polynomials evaluated at
    other Misiurewicz parameters. Goksel \cite{GokselIterates} studies
    the factorization of iterates of PCF unicritical polynomials of
    prime degree over their fields of definition.

    \subsection*{Proof strategy}
    Our reductions to roots of unity use two consequences of
    Kronecker's theorem (Theorem~\ref{thm:Kronecker}).
    Write $\OO_F$ for the ring of integers of a number field $F$.
    \begin{theorem}[Small CM integers]\label{thm:CMgapIntro}
        Let $F$ be a CM field and $\alpha\in\OO_F$. If
        $|\sigma(\alpha)|<\sqrt2$ for every embedding
        $\sigma:F\hookrightarrow\C$, then one of the following alternatives holds
        \begin{equation}\label{eq:small-cm-integer-dichotomy}
            \alpha=0\qquad\text{or}\qquad\alpha\in\mui.
        \end{equation}
        The bound $\sqrt2$ is sharp.
    \end{theorem}

    Section~\ref{sec:integers-high} proves this theorem and uses the
    escape bound $|\sigma(c)|\le2^{1/(d-1)}$ to settle $d\ge3$ by a
    calculation with roots of unity. At $d=3$, Lemma~\ref{lem:cubic}
    excludes equality at $\sqrt2$.

    For $d=2$ the bound is $2$. The fact that a CM unit divided by its
    complex conjugate is a root of unity
    \cite[proof of Theorem~4.12]{Washington} gives the following
    result for the unit group $\OO_F^\times$.
    \begin{theorem}[Small CM units]\label{thm:CMunitIntro}
        If $u\in\OO_F^\times$ for a CM field $F$ and
        $|\sigma(u)|\le2$ for every $\sigma:F\hookrightarrow\C$, then
        $u$ admits a decomposition of the form
        \begin{equation}\label{eq:small-cm-unit-decomposition}
            u=\xi_1+\xi_2\qquad\text{for some }\xi_1,\xi_2\in\mui.
        \end{equation}
    \end{theorem}
    Section~\ref{sec:quadratic} proves Theorem~\ref{thm:CMunitIntro}
    and uses it in the exclusion of larger eventual periods described above.

    Section~\ref{sec:proof} uses the totally real classification for
    real parameters and CM integer rigidity for periodic critical
    points. For positive preperiod and eventual period one or two, the
    relevant parameter polynomials are irreducible and have a real root,
    with one exception, the polynomial whose roots are $\pm i$.
    Since a nonreal element of a CM field cannot have a real conjugate,
    only $\pm i$ remain among the nonreal parameters in these cases.
    No general irreducibility conjecture is needed.
    The CM conjugation lemma is proved in Section~\ref{sec:singular},
    which compares Theorem~\ref{thm:main} and Corollary~\ref{cor:ab}
    with singular moduli.

    \section{Comparison with singular moduli}\label{sec:singular}
    A singular modulus is the $j$-invariant of a CM elliptic curve.
    The analogy of Baker and DeMarco \cite[\S1.2]{BakerDeMarco} suggests
    comparing which singular moduli and which PCF parameters lie in CM
    fields or in $\Q^{\ab}$.

    To use real conjugates, we recall how the canonical involution
    $\iota$ of a CM field $F$ acts under embeddings \cite[p.~39]{Washington}.
    Write $F^+$ for its maximal totally real subfield. The following lemma
    also applies in Sections~\ref{sec:integers-high}--\ref{sec:proof}.
    \begin{lemma}[CM conjugation and real conjugates]\label{lem:CMconj}
        For every complex embedding $\sigma:F\hookrightarrow\C$ one has
        \begin{equation}\label{eq:cm-conjugation-embeddings}
            \sigma(\iota x)=\overline{\sigma(x)}\qquad(x\in F).
        \end{equation}
        If $x\in F$ has a real Galois conjugate, then $x\in F^+$ and
        is totally real.
    \end{lemma}
    \begin{proof}
        Write $F=F^+(\theta)$ with $\iota\theta=-\theta\ne0$. Each real
        embedding of $F^+$ sends $\theta^2$ to a negative number, since
        a positive image would give a real embedding of $F$.
        Thus $\sigma(\theta)$ is purely imaginary, and writing
        $x=a+b\theta$ with $a,b\in F^+$ gives the identity.
        An embedding of $\Q(x)$ with real image extends to some $\sigma$
        on $F$. Then $\sigma(x)=\sigma(\iota x)$ forces $x\in F^+$.
    \end{proof}

    For an imaginary quadratic order $\OO$, let $\Cl(\OO)$ be its group
    of proper invertible ideal classes, with exponent
    $\operatorname{exp}\Cl(\OO)$, the least positive $e$ such that
    $a^e=1$ for all $a\in\Cl(\OO)$. The trivial group has exponent $1$.
    Ring class fields connect this exponent to the field of a singular
    modulus.

    \begin{theorem}[Ring class fields, {\cite[Theorem~2.1(1)--(4)]{DanielsLozano}},
        {\cite[Proposition~2.1]{Voight}}]\label{thm:CMinput}
        Let $\OO$ be an order in an imaginary quadratic field $K$, and let
        $E/\C$ have complex multiplication by $\OO$. Let $x=j(E)$ be its
        $j$-invariant, and put $h=\#\Cl(\OO)$. Then $x$ is an algebraic integer,
        $L=K(x)$ is the ring class field of $\OO$, and the degrees satisfy
        \begin{equation}\label{eq:ring-class-degrees}
            [\Q(x):\Q]=[L:K]=h.
        \end{equation}
        The $j$-invariants of the elliptic curves with complex multiplication
        by $\OO$ form one Galois orbit over $K$.
        The extension $L/\Q$ is Galois, and, writing $\iota$ for complex
        conjugation, one has
        \begin{equation}\label{eq:ring-class-galois-action}
            \begin{aligned}
                \Gal(L/K)&\cong\Cl(\OO),\\
                \Gal(L/\Q)&\cong\Cl(\OO)\rtimes\langle\iota\rangle,\\
                \iota a\iota^{-1}&=a^{-1}\qquad(a\in\Cl(\OO)).
            \end{aligned}
        \end{equation}
    \end{theorem}
    The first reference gives the assertions about the $j$-invariants
    and the degrees. The second describes the Galois group.
    For a singular modulus
    $x=j(E)$ with endomorphism order $\OO$, we obtain the classical equivalences
    \begin{equation}\label{eq:singular}
        \begin{aligned}
            &x\text{ belongs to some CM number field}\\[-2pt]
            &\quad\Longleftrightarrow\quad x\in\Q^{\ab}\\[-2pt]
            &\quad\Longleftrightarrow\quad
            \operatorname{exp}\Cl(\OO)\le2.
        \end{aligned}
    \end{equation}
    The principal class gives a real conjugate $j(\C/\OO)$ of $x$,
    since the lattice $\OO$ is stable under complex conjugation.
    If $x$ lies in a CM field, Lemma~\ref{lem:CMconj} shows that $x$
    is totally real. Put $K=\OO\otimes_\Z\Q$, $L=K(x)$, and $H=\Gal(L/\Q(x))$.
    Theorem~\ref{thm:CMinput} gives $|H|=2$. If $x$ is totally real,
    then each $g^{-1}\iota g$ for $g\in\Gal(L/\Q)$ fixes $x$, hence is
    the unique nonidentity element $\iota$ of $H$. Thus $\iota$ is central
    and its inversion action forces $\operatorname{exp}\Cl(\OO)\le2$. Conversely, that
    exponent condition makes $\Gal(L/\Q)$ abelian, hence $x\in\Q^{\ab}$.
    Every element of $\Q^{\ab}$ lies in a CM field by
    Theorem~\ref{thm:KW} and the observation following it.
    Corollary~\ref{cor:ab} is the analogue of the first equivalence
    in \eqref{eq:singular} for PCF parameters.

    To deduce finiteness, we recall an unconditional result that includes
    nonmaximal orders.
    \begin{theorem}[Finiteness for imaginary quadratic orders, {\cite[Proposition~8.1]{Voight}}]\label{thm:classExponent}
        There are only finitely many imaginary quadratic orders $\OO$, up to
        isomorphism, for which $\operatorname{exp}\Cl(\OO)\le2$.
    \end{theorem}
    Theorems~\ref{thm:CMinput} and~\ref{thm:classExponent} with
    \eqref{eq:singular} give finiteness of singular moduli in CM fields.
    Theorem~\ref{thm:main} gives the analogous explicit PCF parameter
    lists for each $d\ge2$.

    \section{CM integers and the case $d\ge3$}\label{sec:integers-high}
    We prove Theorem~\ref{thm:CMgapIntro} and combine it with critical
    orbit bounds to settle the case $d\ge3$ of Theorem~\ref{thm:main}.

    \subsection{CM integers and critical orbit bounds}\label{sec:general}
    \begin{theorem}[Kronecker, {\cite[I, p.~173]{Kronecker}}]\label{thm:Kronecker}
        Let $\alpha$ be an algebraic integer. If every complex conjugate of
        $\alpha$ has absolute value $1$, then $\alpha$ is a root of unity.
    \end{theorem}
    \begin{proof}[Proof of Theorem~\ref{thm:CMgapIntro}]
        Let $F$ be CM and $\alpha\in\OO_F$ satisfy
        $|\sigma(\alpha)|<\sqrt2$ for every embedding $\sigma:F\hookrightarrow\C$.
        The case $\alpha=0$ is immediate. Otherwise put
        $\beta=\alpha\iota(\alpha)\in\OO_{F^+}$.
        Each $\tau:F^+\hookrightarrow\R$ extends to an embedding $\sigma$
        on $F$, so Lemma~\ref{lem:CMconj} gives
        \begin{equation}\label{eq:cm-gap-positive-conjugates}
            0<\tau(\beta)=|\sigma(\alpha)|^2<2.
        \end{equation}
        If $\beta\ne1$, its difference from $1$ has nonzero integer norm, yet
        \begin{equation}\label{eq:cm-gap-norm-contradiction}
            0<\left|\Nm_{F^+/\Q}(\beta-1)\right|
            =\prod_{\tau:F^+\hookrightarrow\R}|\tau(\beta)-1|<1.
        \end{equation}
        Thus $\beta=1$, so all conjugates of $\alpha$ have absolute value
        $1$ and Theorem~\ref{thm:Kronecker} gives $\alpha\in\mui$.
        The example $1+i\in\Q(i)$ has both conjugates of absolute value
        $\sqrt2$ and is not a root of unity, proving optimality.
    \end{proof}

    The degree-dependent bound below also follows from Buff's estimates
    for $w\mapsto aw^d+1$ \cite[Lemma~7 and Corollary~8]{BuffUnicritical}.
    For $c\ne0$, the change of coordinates $z=cw$ gives $a=c^{d-1}$,
    so the bound $|aw^{d-1}|\le2$ along the critical orbit becomes
    $|z|^{d-1}\le2$. The coarser bound $2$ in arbitrary degree is given in
    \cite[Lemma~4.1]{FerragutiPagano}. We give the integrality and escape
    argument in the coordinates of $f_{d,c}$.
    \begin{lemma}[Uniform escape bound]\label{lem:escape}
        Let $d\ge2$ and suppose $f_{d,c}$ is PCF. Then $c$ is an
        algebraic integer. For every embedding $\sigma:\Q(c)\hookrightarrow\C$
        and every $r\ge1$, the critical iterates satisfy the bound
        \begin{equation}\label{eq:critical-orbit-bound}
            \bigl|\sigma(f_{d,c}^{\circ r}(0))\bigr|\le 2^{1/(d-1)}.
        \end{equation}
        All conjugates of $c$ have the same exact critical orbit type.
    \end{lemma}
    \begin{proof}
        If $(m,n)$ is the exact critical orbit type, then
        $f_{d,c}^{\circ(m+n)}(0)-f_{d,c}^{\circ m}(0)=0$ is a monic
        polynomial equation for $c$, so $c$ is an algebraic integer.
        Embeddings preserve equality and distinctness of orbit points, so it
        suffices to bound one orbit.

        If $|c|>2^{1/(d-1)}$, then every $z$ with $|z|=s\ge|c|$
        satisfies the estimate
        \begin{equation}\label{eq:escape-large-parameter}
            |f_{d,c}(z)|\ge s^d-|c|\ge s(s^{d-1}-1)>s.
        \end{equation}
        Starting at $c$, the absolute values would therefore strictly
        increase. If $|c|\le2^{1/(d-1)}$ and an orbit point $z$ has
        $|z|=s>2^{1/(d-1)}$, then
        \begin{equation}\label{eq:escape-large-orbit-value}
            |f_{d,c}(z)|\ge s^d-|c|>s,
        \end{equation}
        because $s^d-s$ is strictly increasing for $s\ge2^{1/(d-1)}$
        and its endpoint value satisfies
        \begin{equation}\label{eq:escape-threshold}
            2^{d/(d-1)}-2^{1/(d-1)}=2^{1/(d-1)}\ge|c|.
        \end{equation}
        Thus in either case the orbit has strictly increasing absolute
        values from that point on, contradicting finiteness.
    \end{proof}

    For $d\ge4$, Lemma~\ref{lem:escape} gives the strict inequality
    required in Theorem~\ref{thm:CMgapIntro}. In degree three we must
    exclude equality.
    \begin{lemma}[Strict bound for cubics]\label{lem:cubic}
        If $f_{3,c}$ is PCF, then every conjugate $c'$ of $c$ satisfies
        $|c'|<\sqrt2$.
    \end{lemma}
    \begin{proof}
        Lemma~\ref{lem:escape} bounds $|c'|$ by $\sqrt2$.
        Suppose equality holds. The triangle inequality then gives
        \begin{equation}\label{eq:cubic-boundary-second-iterate}
            |f_{3,c'}^{\circ 2}(0)|=|c'^3+c'|
            \ge |c'|^3-|c'|=\sqrt2.
        \end{equation}
        A strict inequality would contradict Lemma~\ref{lem:escape}.
        Equality forces $(c')^2=-2$,
        so $f_{3,c'}^{\circ 2}(0)=-c'$ and $f_{3,c'}^{\circ 3}(0)=3c'$.
        The last equality contradicts Lemma~\ref{lem:escape}, since
        $|3c'|=3\sqrt2$.
    \end{proof}

    \subsection{Classification for $d\ge3$}\label{sec:high}
    \begin{proposition}\label{prop:high}
        Let $d\ge3$ and suppose $c$ lies in a CM field. Then $f_{d,c}$
        is PCF if and only if $c\in\mathcal{S}_d$ as defined in
        \eqref{eq:Sd}. The exact critical orbit types are $(0,1)$ for $c=0$,
        $(0,2)$ for $c^{d-1}=-1$, $(2,1)$ for
        $c^{d-1}\in\{\zeta_3,\zeta_3^2\}$
        when $d\equiv0\pmod6$, and $(0,3)$ for the latter parameters
        when $d\equiv1\pmod6$.
    \end{proposition}
    \begin{proof}
        By Lemmas~\ref{lem:escape} and \ref{lem:cubic}, every Galois
        conjugate of a PCF parameter $c$ has absolute value strictly less than
        $\sqrt2$ for $d\ge3$.  Theorem~\ref{thm:CMgapIntro}
        therefore gives $c=0$ or $c\in\mui$.
        The case $c=0$ is clear. Suppose $c\in\mui$ and put
        $t=c^{d-1}\in\mui$. The first critical iterates are
        \begin{equation}\label{eq:root-unity-first-iterates}
            \begin{aligned}
                f_{d,c}(0)&=c,\\
                f_{d,c}^{\circ 2}(0)&=c(1+t),\\
                f_{d,c}^{\circ 3}(0)&=c\left(1+(1+t)^dt\right).
            \end{aligned}
        \end{equation}
        If $t$ has order $N\ge4$, replace $c$ by the conjugate obtained
        by extending the embedding $t\mapsto\zeta_N$ of $\Q(t)$ to
        $\Q(c)$. The new parameter is still PCF by Lemma~\ref{lem:escape},
        and its second iterate satisfies
        \begin{equation}\label{eq:high-order-second-iterate}
            |f_{d,c}^{\circ 2}(0)|=|1+t|=2\cos(\pi/N)\ge\sqrt2.
        \end{equation}
        As $|c|=1$ and $d\ge3$, we also obtain the inequalities
        \begin{equation}\label{eq:high-order-third-iterate}
            |f_{d,c}^{\circ 3}(0)|\ge|f_{d,c}^{\circ 2}(0)|^d-1
            \ge2\sqrt2-1>\sqrt2\ge2^{1/(d-1)}.
        \end{equation}
        This contradicts Lemma~\ref{lem:escape}. Hence the order of
        $t$ is at most three. If $t=1$, the second critical iterate equals
        $2c$ and the orbit escapes. If $t=-1$, the second critical iterate is $0$,
        so the critical point has exact period two.

        It remains to consider $t\in\{\zeta_3,\zeta_3^2\}$.
        Since $1+t=-t^2$, we have
        \begin{equation}\label{eq:cubicroot}
            \begin{aligned}
                f_{d,c}^{\circ 2}(0)&=-ct^2,\\
                f_{d,c}^{\circ 3}(0)&=c\bigl(1+(-1)^d t^{2d+1}\bigr).
            \end{aligned}
        \end{equation}
        For odd $d$, formula \eqref{eq:cubicroot} yields
        $f_{d,c}^{\circ 3}(0)=0$ if $d\equiv1\pmod6$,
        in which case the critical point has exact period three.
        For $d\equiv3,5\pmod6$, we have the estimate
        $|f_{d,c}^{\circ 3}(0)|=|1-t^{2d+1}|=\sqrt3>2^{1/(d-1)}$,
        so the orbit is not finite.

        For even $d$, formula \eqref{eq:cubicroot} gives
        $f_{d,c}^{\circ 3}(0)=f_{d,c}^{\circ 2}(0)$ if $d\equiv0\pmod6$,
        in which case the exact critical orbit type is $(2,1)$.
        If $d\equiv4\pmod6$, it gives $f_{d,c}^{\circ 3}(0)=2c$.
        If $d\equiv2\pmod6$, then we obtain
        $f_{d,c}^{\circ 3}(0)=-ct$ and
        $f_{d,c}^{\circ 4}(0)=c(1+t^{d+1})=2c$.
        In both cases the absolute value $2>2^{1/(d-1)}$ contradicts
        Lemma~\ref{lem:escape}. The orbit computations also verify
        that each listed parameter is PCF.
    \end{proof}

    \begin{remark}[Parameter counts and conjugacy classes]\label{rem:high}
        The solutions of $c^{d-1}=-1$ are $(2d-2)$th roots of unity.
        For $d\equiv0,1\pmod6$ the two further equations
        $c^{d-1}=\zeta_3,\zeta_3^2$ contribute $2(d-1)$ distinct
        parameters, yielding \eqref{eq:cardinality}.
        In the moduli space of unicritical polynomials, the scaling
        $c\mapsto\xi c$ with $\xi^{d-1}=1$ identifies the
        parameters having the same $t=c^{d-1}$. For $d\ge3$, the lists
        therefore give four affine conjugacy classes when $d\equiv0,1\pmod6$
        and two otherwise. For $d\equiv1\pmod6$, the period-three cases above
        correspond to the factor $a^2+a+1$ of Buff's period-three Gleason
        polynomial in the parameter $a=c^{d-1}$ \cite[Proposition~5]{BuffUnicritical}.
    \end{remark}
    \section{The quadratic case and roots of unity}\label{sec:quadratic}
    Write $f_c=f_{2,c}$. We rule out parameters in CM fields whose
    critical orbit is strictly preperiodic with eventual period at least
    three. The proof reduces to four cyclotomic families.
    Section~\ref{sec:proof} treats the remaining cases.

    \subsection{Unit theorem and reduction to roots of unity}\label{sec:prelim}\label{sec:CMrigidity}
    Suppose $f_c$ has exact critical orbit type $(m,n)$.
    If $m=1$, then $f_c^{\circ(n+1)}(0)=c$ forces $f_c^{\circ n}(0)=0$,
    since $f_c(z)=c$ only at $z=0$. This would make $0$ periodic,
    so a positive preperiod must be at least two.
    To apply Theorem~\ref{thm:CMunitIntro} to $c$ and $c+1$, we
    first bound their conjugates.

    \begin{lemma}[Archimedean bounds]\label{lem:arch}
        Let $c$ be a PCF parameter for $f_{2,c}$. For every embedding
        $\sigma:\Q(c)\hookrightarrow\C$ and every $r\ge1$, we have the inequalities
        \begin{equation}\label{eq:quadratic-orbit-bounds}
            |\sigma(c)|\le2,\qquad |\sigma(f_c^{\circ r}(0))|\le2.
        \end{equation}
        If $c\ne-2$, all these inequalities are strict. In all cases,
        every such embedding $\sigma$ satisfies $|\sigma(c+1)|<2$.
    \end{lemma}
    \begin{proof}
        Lemma~\ref{lem:escape} gives the weak bounds for every conjugate.
        It suffices to check strictness for $c$ itself, since its conjugates
        are PCF and none equals $-2$ unless $c=-2$.
        If $|c|=2$, then $|f_c^{\circ2}(0)|=|c(c+1)|\ge2$, with equality
        only at $c=-2$. Thus $|c|<2$ for $c\ne-2$. In that case, if an
        iterate has absolute value $2$, the next has absolute value at
        least $4-|c|>2$. For $|c|<1$, the triangle inequality gives
        $|c+1|<2$. If $|c|\ge1$ and $c\ne-2$, we have
        \begin{equation}\label{eq:quadratic-shift-bound}
            |c+1|=\frac{|f_c^{\circ 2}(0)|}{|c|}<2.
        \end{equation}
        For $c=-2$, one has $|c+1|=1$.
    \end{proof}

    Goksel's theorem supplies the unit condition. We include a proof
    for convenience.
    \begin{theorem}[Goksel, {\cite[Theorem~1.3(a)]{GokselOrbit}}]\label{thm:Goksel}
        Let $d\ge2$ and suppose $f_{d,c}$ has exact critical orbit type
        $(m,n)$ with $m\ge1$. For every integer $r\ge1$ with $n\nmid r$,
        one has
        \begin{equation}\label{eq:critical-orbit-units}
            f_{d,c}^{\circ r}(0)\in\OO_{\Q(c)}^\times.
        \end{equation}
        The proof below also shows that these elements are units in $\Z[c]$.
        In particular, if $d=2$ and $n\ge3$, then both $c$ and $c+1$ are units.
    \end{theorem}
    \begin{proof}
        By Lemma~\ref{lem:escape}, $A=\Z[c]\subseteq\OO_{\Q(c)}$.
        Put $a=f_{d,c}^{\circ r}(0)$.
        As $f_{d,c}^{\circ r}(z)-a\in zA[z]$, induction gives
        $f_{d,c}^{\circ(kr)}(0)\in aA$ for $k\ge1$.
        Choose a positive multiple $\ell\ge m$ of $r$ and put
        $x_j=f_{d,c}^{\circ(\ell+j)}(0)$ for $j\ge0$. These points have
        exact period $n$, so $\Delta_j=x_{j+r}-x_j\ne0$ since $n\nmid r$.
        Factoring the difference of $d$th powers gives
        \begin{equation}\label{eq:orbit-difference-recurrence}
            \Delta_{j+1}=D_j\Delta_j,
            \qquad
            D_j=\sum_{h=0}^{d-1}x_{j+r}^{\,d-1-h}x_j^h\in A.
        \end{equation}
        Multiplying over $j=0,\ldots,n-1$ and using $\Delta_n=\Delta_0$,
        we cancel the nonzero product of the $\Delta_j$ in $A$ to obtain
        \begin{equation}\label{eq:orbit-difference-product}
            \prod_{j=0}^{n-1}D_j=1.
        \end{equation}
        Thus $D_0$ is a unit of $A$. Both $x_0$ and $x_r$ lie in $aA$,
        so $D_0\in a^{d-1}A\subseteq aA$ and $a$ is a unit.
        For $d=2$ and $n\ge3$, take $r=2$ to see that each factor of
        the unit $f_c^{\circ2}(0)=c(c+1)$ is a unit in $A$.
    \end{proof}

    \begin{proof}[Proof of Theorem~\ref{thm:CMunitIntro}]
        Let $F$ be a CM field and $u\in\OO_F^\times$ satisfy
        $|\sigma(u)|\le2$ for every embedding $\sigma:F\hookrightarrow\C$.
        The quotient $u/\iota u$ is an algebraic unit whose conjugates
        have absolute value $1$ by Lemma~\ref{lem:CMconj}.
        Thus $u/\iota u=\eta\in\mui$ by Theorem~\ref{thm:Kronecker}, as in
        \cite[proof of Theorem~4.12]{Washington}.
        Choose $\varepsilon\in\mui$ with $\varepsilon^2=\eta$ and put
        $v=u/\varepsilon$, which is integral because $\varepsilon$ is a unit.
        For every embedding $\sigma:F(\varepsilon)\hookrightarrow\C$,
        Lemma~\ref{lem:CMconj} and $|\sigma(\varepsilon)|=1$ give
        \begin{equation}\label{eq:unit-twist-real}
            \overline{\sigma(v)}=\sigma(\iota u)\sigma(\varepsilon)
            =\frac{\sigma(u)}{\sigma(\varepsilon)}=\sigma(v).
        \end{equation}
        Thus every conjugate of $v$ is real and has absolute value
        $|\sigma(u)|\le2$. A root $w$ of $T^2-vT+1$ is integral, and all
        its conjugates lie on the unit circle since the conjugates of $v$
        lie in $[-2,2]$. By Theorem~\ref{thm:Kronecker}, $w\in\mui$.
        Since $v=w+w^{-1}$, we obtain
        \begin{equation}\label{eq:small-unit-sum}
            u=\varepsilon(w+w^{-1})
            =\varepsilon w+\varepsilon w^{-1}.
        \end{equation}
    \end{proof}

    Applying the decomposition to $c$ and $c+1$ turns $(c+1)-c=1$
    into a five-term vanishing sum. Such a sum is \emph{minimal}
    if no nonempty proper subcollection of terms sums to zero,
    counting terms with multiplicity.
    \begin{theorem}[Lam--Leung, short-sum consequence of {\cite[Theorem~4.8]{LamLeung}}]\label{thm:shortsums}
        Suppose $\alpha_1+\cdots+\alpha_\ell=0$ is a minimal vanishing sum
        of roots of unity with $1\le\ell\le5$. Then $\ell$ belongs to
        $\{2,3,5\}$. After reordering the terms, there is a root of unity
        $\varepsilon$ for which the terms have the representation
        \begin{equation}\label{eq:minimal-root-sum}
            \alpha_j=\varepsilon\zeta_\ell^{\,j-1}
            \qquad(1\le j\le\ell).
        \end{equation}
    \end{theorem}
    \begin{proof}
        The cited theorem says that a minimal vanishing
        sum is either a rotated regular polygon of prime order, or the
        least common multiple of the orders of its terms has at least three
        distinct prime divisors and the sum has length at least
        $p_1(p_2-1)+p_3-p_2\ge6$, where $p_1<p_2<p_3$ are the three
        smallest such primes. Thus the latter alternative is
        impossible here.
    \end{proof}
    \begin{lemma}[Reduction to four families]\label{lem:five}
        Let $c$ lie in a CM field and suppose $f_c$ has exact critical
        orbit type $(m,n)$ with $m\ge2$ and $n\ge3$.
        Then $c$ satisfies $c(c+1)\ne0$ and belongs to one of the following four families
        \begin{equation}\label{eq:fourfamilies}
            c\in\mui,\qquad
            c=\zeta-1,\qquad
            c=\zeta+\zeta_3^e,\qquad
            c=\zeta_5^a+\zeta_5^b\quad(1\le a<b\le4),
        \end{equation}
        where $\zeta\in\mui$ and $e\in\{1,2\}$.
    \end{lemma}
    \begin{proof}
        By Theorem~\ref{thm:Goksel}, $c$ and $c+1$ are units in
        $\OO_{\Q(c)}$ and therefore in the ring of integers of any CM
        field containing $c$.
        By Lemma~\ref{lem:arch}, their conjugates have absolute
        value at most $2$. Thus Theorem~\ref{thm:CMunitIntro} gives
        $c=r+s$ and $c+1=t+u$ with $r,s,t,u\in\mui$, so
        \begin{equation}\label{eq:five-root-relation}
            1+r+s-t-u=0.
        \end{equation}
        All terms are roots of unity. Decompose the sum into minimal
        vanishing subsums. By Theorem~\ref{thm:shortsums}, these form a
        pentagon or an opposite pair and a triangle. In the pentagon case, $1$
        fixes the rotation, so
        $r,s$ are distinct nontrivial fifth roots. In the pair and triangle
        case, pairing $1$ with $r$ or $s$ gives $c=\zeta-1$.
        Pairing it with $-t$ or $-u$ gives $c\in\mui$.
        The pairs $(r,s)$ and $(-t,-u)$ would give $c=0$ and $c=-1$,
        respectively, excluded by the unit conditions. The only other
        pairing, up to symmetry, gives $r=t$ and $1+s-u=0$.
        Since $|s|=|1+s|=1$, we get $s\in\{\zeta_3,\zeta_3^2\}$.
    \end{proof}

    \subsection{Exclusion of the cyclotomic families}\label{sec:cyclotomic}
    We analyze the four families in \eqref{eq:fourfamilies} to prove
    Proposition~\ref{prop:nohigher}, beginning with roots of unity.
    \begin{lemma}[Root-of-unity parameters]\label{lem:rootunity}
        If $c$ is a root of unity and $f_c$ is PCF, then
        $c$ belongs to $\{-1,i,-i\}$.
    \end{lemma}
    \begin{proof}
        Put $\gamma=(1+\sqrt5)/2$.
        If $c$ has order $N\ge6$, replace it by its conjugate $\zeta_N$,
        which is also a PCF parameter by Lemma~\ref{lem:escape}. Its second
        iterate satisfies
        \begin{equation}\label{eq:root-unity-escape-threshold}
            |f_c^{\circ 2}(0)|=|c+1|=2\cos(\pi/N)>
            \gamma.
        \end{equation}
        Since $|c|=1$ and $r^2-1>r$ for $r>\gamma$, subsequent absolute
        values strictly increase, contradicting finiteness of the orbit.
        For $N=5$, pass to the conjugate $c=\zeta_5$, which satisfies
        $f_c^{\circ 3}(0)=c^3-1$.
        The equality $|f_c^{\circ 3}(0)|^2=(5+\sqrt5)/2$ gives the further bound
        $|f_c^{\circ 4}(0)|\ge|f_c^{\circ 3}(0)|^2-1=(3+\sqrt5)/2>2$.
        For a primitive cubic root, the critical iterates satisfy $f_c^{\circ 4}(0)=2c$
        and $|f_c^{\circ 5}(0)|^2=13$.
        The parameter $c=1$ also has an escaping critical orbit.
        The only remaining orders are $2$ and $4$, giving $-1,i,-i$.
        Direct iteration gives exact critical orbit types $(0,2)$ for
        $c=-1$ and $(2,2)$ for $c=\pm i$.
    \end{proof}

    For the next two families, we average the squared absolute values of
    a critical iterate over its conjugates. If the mean exceeds $4$, some
    conjugate has absolute value greater than $2$, contradicting
    Lemma~\ref{lem:arch}. For $u$ in a number field $L$, set
    \begin{equation}\label{eq:embedding-mean-definition}
        \av|u|^2=\frac{1}{[L:\Q]}
        \sum_{\sigma:L\hookrightarrow\C}|\sigma(u)|^2.
    \end{equation}
    The mean is independent of $L$, since each embedding of $\Q(u)$
    has $[L:\Q(u)]$ extensions to $L$.

    Write $\mu$ for the M\"obius function and $\varphi$ for Euler's totient.
    To compute the averages, we use the normalized Ramanujan sum
    \begin{equation}\label{eq:ramanujan-mean-definition}
        A_N(j)=\frac{1}{\varphi(N)}
        \sum_{\substack{1\le k\le N\\\gcd(k,N)=1}}\zeta_N^{kj}
        \qquad(j\in\Z).
    \end{equation}
    With the convention $\mu(1)=\varphi(1)=1$, the classical
    Ramanujan-sum identity takes the following form.
    \begin{lemma}[Cyclotomic mean formula]\label{lem:means}
        For $q=N/\gcd(N,j)$, the mean satisfies the identity
        \begin{equation}\label{eq:ramanujan}
            A_N(j)=\frac{\mu(q)}{\varphi(q)}.
        \end{equation}
    \end{lemma}
    \begin{proof}
        The primitive $N$th roots map uniformly onto the primitive $q$th
        roots under $z\mapsto z^j$, since reduction
        $(\Z/N\Z)^\times\to(\Z/q\Z)^\times$ is surjective, with all fibers
        of size $\varphi(N)/\varphi(q)$. Their sum is $\mu(q)$, by M\"obius
        inversion. Dividing by $\varphi(q)$ gives the formula.
    \end{proof}

    \begin{lemma}[The family $\zeta-1$]\label{lem:zetaminus}
        If $\zeta$ is a root of unity of order $N\notin\{1,2,6\}$, then
        $c=\zeta-1$ has a conjugate $c'$ with $|f_{c'}^{\circ 3}(0)|>2$.
    \end{lemma}
    \begin{proof}
        One computes
        \begin{equation}\label{eq:zeta-minus-third-iterate}
            f_{\zeta-1}^{\circ 3}(0)=\zeta^4-2\zeta^3+\zeta^2+\zeta-1.
        \end{equation}
        Averaging the squared absolute value gives
        \begin{equation}\label{eq:Tminus}
            T_N=\av|f_{\zeta-1}^{\circ 3}(0)|^2
            =8-8A_N(1)-4A_N(2)+6A_N(3)-2A_N(4).
        \end{equation}
        Write $N=2^a3^bR$ with $a,b\ge0$, $R\ge1$, and $\gcd(R,6)=1$.
        If $R$ is not squarefree, all four means vanish and $T_N=8$.
        Otherwise let $h=\mu(R)/\varphi(R)$. Formula~\eqref{eq:ramanujan}
        and its multiplicativity for coprime moduli give
        \begin{equation}\label{eq:ramanujan-factorization}
            A_N(j)=A_{2^a}(j)A_{3^b}(j)h\qquad(1\le j\le4),
        \end{equation}
        since $\gcd(R,j)=1$ in this range. For example, when $(a,b)=(0,1)$,
        the means $A_N(1),A_N(2),A_N(4)$ all equal $-h/2$, while $A_N(3)=h$.
        Substituting in \eqref{eq:Tminus} gives $T_N=8+13h$ in this case.
        The same substitution gives the remaining rows of the table.
        \begin{equation}\label{eq:zeta-minus-mean-table}
            \begin{array}{c@{\qquad}l}
                \toprule
                (a,b)&T_N\\\midrule
                (0,0)&8-8h\\
                (1,0)&8-4h\\
                (0,1)&8+13h\\
                (1,1)&8-7h\\
                (2,0),(3,0)&8+2h\\
                (2,1),(3,1)&8-h\\
                (0,2)&8-3h\\
                (1,2)&8+3h\\
                \text{all remaining pairs}&8\\\bottomrule
            \end{array}
        \end{equation}
        If $a\ge2$ or $b\ge2$, this gives $T_N\ge5$.
        For the first four rows with $R>1$, every prime factor of $R$ is at
        least $5$, so the estimates $|h|\le1/4$ and $T_N\ge19/4>4$ follow.
        For $R=1$, these four rows correspond to $N=1,2,3,6$, with means
        $0,4,21,1$, respectively. The bound $T_N>4$ therefore holds for
        every $N\notin\{1,2,6\}$.
    \end{proof}

    \begin{lemma}[The family $\zeta+\zeta_3^e$]\label{lem:zetaplus}
        Let $\zeta\in\mui$ and $e\in\{1,2\}$, and put
        $c=\zeta+\zeta_3^e$. If $f_{2,c}$ is PCF, then $c\in\{0,-1\}$.
    \end{lemma}
    \begin{proof}
        Let $N$ be the order of $\zeta$. Passing to a conjugate preserves
        the PCF condition and fixes $0,-1$, so we may assume $\zeta=\zeta_N$
        after adjusting $e$ within $\{1,2\}$.
        Put $B_e=1+2\zeta_3^e$,
        so $B_e^2=-3$ and $\overline{B_e}=-B_e$. Then
        \begin{equation}\label{eq:zeta-plus-second-iterate}
            \begin{aligned}
                f_c^{\circ 2}(0)&=\zeta^2+B_e\zeta-1,\\
                |f_c^{\circ 2}(0)|^2&=5-\zeta^2-\zeta^{-2}
                -2B_e(\zeta-\zeta^{-1}).
            \end{aligned}
        \end{equation}
        First suppose $3\nmid N$. The images of $\zeta_N$ and $\zeta_3$
        vary independently under Galois automorphisms, so the mean is
        \begin{equation}\label{eq:zeta-plus-coprime-mean}
            T_N=\av|f_c^{\circ 2}(0)|^2=5-2A_N(2).
        \end{equation}
        For $N\ge3$, the integer $q=N/\gcd(N,2)$ satisfies $q>1$ and $3\nmid q$.
        Formula~\eqref{eq:ramanujan} gives $T_N=5-2\mu(q)/\varphi(q)$.
        For $\mu(q)\le0$, the estimate $T_N\ge5>4$ is immediate.
        If $\mu(q)>0$, then $q$ is squarefree with an even number of
        distinct prime factors. As $q>1$ and $3\nmid q$, it has a prime
        divisor $p\ge5$. The bounds $\varphi(q)\ge p-1\ge4$ and
        $T_N\ge9/2>4$ follow.
        For the remaining orders $N=1$ and $2$, direct calculation gives
        $|f_{1+\zeta_3^e}^{\circ 3}(0)|^2=7$ and
        $|f_{-1+\zeta_3^e}^{\circ 3}(0)|^2=21$, respectively.

        Now suppose $3\mid N$. We then have $\zeta_3^e=\zeta_N^{eN/3}$.
        The relation $1+\zeta_3^e+\zeta_3^{2e}=0$ gives
        \begin{equation}\label{eq:Te}
            T_{N,e}=\av|f_c^{\circ 2}(0)|^2
            =5-2A_N(2)-4A_N(1)-8A_N(1+eN/3).
        \end{equation}
        If $9\mid N$, the three Ramanujan means vanish and the mean
        satisfies $T_{N,e}=5>4$.
        Otherwise write $N=3M$, $3\nmid M$, and put
        \begin{equation}\label{eq:zeta-plus-auxiliary-means}
            h=\frac{\mu(M)}{\varphi(M)},\qquad
            b=\frac{\mu(M/\gcd(M,2))}{\varphi(M/\gcd(M,2))}.
        \end{equation}
        Using \eqref{eq:ramanujan} in \eqref{eq:Te}, we find
        \begin{equation}\label{eq:Tsplit}
            T_{N,e}=\begin{cases}
                5+b+6h,&eM\equiv1\pmod3,\\
                5+b-6h,&eM\equiv2\pmod3.
            \end{cases}
        \end{equation}
        If $M$ is not squarefree, then $h=0$ and $T_{N,e}=5+b\ge4$.
        Equality requires $b=-1$, hence $M/\gcd(M,2)=2$ and $M=4$.
        Now suppose $M$ is squarefree. For $\varphi(M)>7$, the relation
        $b=\pm h$ implies $T_{N,e}\ge5-7/\varphi(M)>4$.
        For squarefree $M$ with $3\nmid M$ and $\varphi(M)\le7$, every odd
        prime divisor is $5$ or $7$, and both cannot occur together.
        Thus the remaining possibilities are $M\in\{1,2,5,7,10,14\}$.
        Together with $M=4$, these give the exact means
        \begin{equation}\label{eq:zeta-plus-exceptional-means}
            \begin{array}{c|cc}
                \toprule
                N&T_{N,1}&T_{N,2}\\\midrule
                3,6&12&0\\12&4&4\\
                15,30&25/4&13/4\\21,42&23/6&35/6\\\bottomrule
            \end{array}
        \end{equation}
        The zero means give $c=-1$ for $(N,e)=(3,2)$ and $c=0$ for
        $(6,2)$. For the six pairs with $0<T_{N,e}\le4$, we test
        individual conjugates by choosing the exponents
        \begin{equation}\label{eq:zeta-plus-conjugating-exponents}
            \begin{array}{c|rrrrrr}
                N&12&12&15&21&30&42\\
                e&1&2&2&1&2&1\\
                k&1&5&2&1&11&1
            \end{array}
        \end{equation}
        Each $k$ defines an automorphism $\zeta_N\mapsto\zeta_N^k$ since
        $\gcd(k,N)=1$. As $ek\equiv1\pmod3$, it sends $c$ to the conjugate
        $c'=e^{i\theta}+\zeta_3$ with $\theta=2\pi k/N$. Formula
        \eqref{eq:zeta-plus-second-iterate} gives
        \begin{equation}\label{eq:theta}
            |f_{c'}^{\circ 2}(0)|^2=3+4\sqrt3\sin\theta+4\sin^2\theta.
        \end{equation}
        Each angle lies in $[\pi/21,5\pi/6]$. Concavity gives
        $\sin\theta\ge3\theta/\pi$ on $[0,\pi/6]$, and
        $\sin\theta\ge1/2$ on $[\pi/6,5\pi/6]$. Thus $\sin\theta\ge1/7$
        for the chosen angles and \eqref{eq:theta} yields
        \begin{equation}\label{eq:zeta-plus-second-iterate-lower-bound}
            |f_{c'}^{\circ 2}(0)|^2\ge3+\frac{4\sqrt3}{7}+\frac4{49}
            >3+\frac{20}{21}+\frac4{49}>4,
        \end{equation}
        contrary to Lemma~\ref{lem:arch}.
    \end{proof}

    \begin{proposition}[Exclusion of the four families]\label{prop:nohigher}
        No parameter $c$ in \eqref{eq:fourfamilies} with $c(c+1)\ne0$
        defines a PCF map $f_{2,c}$ whose critical orbit has eventual
        period $n\ge3$.
    \end{proposition}
    \begin{proof}
        For a root of unity, Lemma~\ref{lem:rootunity} leaves only
        $-1,\pm i$, each of eventual period $2$.
        If $c=\zeta-1$, Lemma~\ref{lem:zetaminus} leaves
        only orders $1,2,6$ for $\zeta$. The corresponding parameters are $0,-2$
        (both of eventual period $1$) and primitive cubic roots
        (which are excluded by Lemma~\ref{lem:rootunity}).
        If $c=\zeta+\zeta_3^e$, Lemma~\ref{lem:zetaplus} leaves $0,-1$,
        contradicting $c(c+1)\ne0$.
        Finally, suppose $c=\zeta_5^a+\zeta_5^b$.
        If the two roots are inverse, one of the conjugates of $c$ equals
        $2\cos(2\pi/5)>0$. For a positive real parameter the critical
        orbit is strictly increasing and hence infinite.
        The four other unordered pairs form one Galois orbit,
        represented by $c=\zeta_5+\zeta_5^2$, for which we have
        \begin{equation}\label{eq:fifth-root-escape}
            |c|=|c+1|=\frac{1+\sqrt5}{2},
            \qquad |f_c^{\circ 2}(0)|=\frac{3+\sqrt5}{2}>2.
        \end{equation}
        This contradicts Lemma~\ref{lem:arch}.
    \end{proof}

    \begin{lemma}[Exclusion of larger eventual periods]\label{lem:excludeHigher}
        There is no parameter $c$ in a CM field for which $f_{2,c}$ has
        exact critical orbit type $(m,n)$ with $m\ge2$ and $n\ge3$.
    \end{lemma}
    \begin{proof}
        Lemma~\ref{lem:five} gives $c(c+1)\ne0$ and places $c$ in
        \eqref{eq:fourfamilies}. Proposition~\ref{prop:nohigher} then
        contradicts $n\ge3$.
    \end{proof}

    \section{Completion of the classification}\label{sec:proof}
    For $f_{2,c}$ over CM fields, Subsection~\ref{sec:prelim} and
    Lemma~\ref{lem:excludeHigher} leave only eventual periods one and two
    when the critical point is strictly preperiodic. If $c$ is real, it is
    totally real by Lemma~\ref{lem:CMconj}, so the following classification
    applies.

    \begin{lemma}[Totally real parameters, {\cite{NP}},
        {\cite[Proposition~1]{BK}}]\label{lem:totreal}
        If $c$ is totally real and $f_c$ is PCF, then
        $c$ belongs to $\{0,-1,-2\}$.
    \end{lemma}
    \begin{proof}
        We recall the elementary argument given by Buff and Koch.
        Lemmas~\ref{lem:escape} and~\ref{lem:arch} show that $c$ is
        integral and every conjugate $c'$ is a PCF parameter in $[-2,2]$.
        Positive $c'$ give strictly increasing orbits, so in fact
        $c'\in[-2,0]$. A root $w$ of $T^2-cT+1$ is integral and all
        its conjugates lie on the unit circle. By Theorem~\ref{thm:Kronecker},
        $w$ is a root of unity, say of order $N$. From $c=w+w^{-1}$,
        we get $2\operatorname{Re}(\zeta_N)\in[-2,0]$, excluding
        $N=1$ and $N\ge5$. The orders $2,3,4$ give $c=-2,-1,0$.
    \end{proof}

    To determine which nonreal parameters with positive preperiod lie
    in CM fields, we look for real roots of their irreducible parameter
    polynomials. For eventual period two, the formula below gives the
    degrees and values needed for this test. We treat eventual period
    one in the proof of Theorem~\ref{thm:main}.
    \begin{theorem}[Goksel {\cite[Corollary~1.5]{Goksel}} and
        Buff--Epstein--Koch {\cite{BEK}}]\label{thm:irred}
        For an indeterminate $c$ and an integer $m\ge2$, put
        $H_m(c)=f_c^{\circ m}(0)-f_c^{\circ(m-1)}(0)+1$ and define
        \begin{equation}\label{eq:Gm2}
            G_{m,2}(c)=
            \begin{cases}
                H_m(c),&m\text{ even},\\
                H_m(c)/(c+1),&m\text{ odd}.
            \end{cases}
        \end{equation}
        Then $G_{m,2}$ is a monic polynomial in $\Z[c]$ that is irreducible
        over $\Q$. Its roots are exactly the parameters for which $f_c$
        has exact critical orbit type $(m,2)$.
    \end{theorem}
    \begin{proof}
        The period-two orbit relation factors as
        \begin{equation}\label{eq:period-two-factorization}
            f_c^{\circ2}(z)-z=(f_c(z)-z)\Phi_2(c,z),
            \qquad \Phi_2(c,z)=z^2+z+c+1.
        \end{equation}
        Direct multiplication gives
        \begin{equation}\label{eq:period-two-transport}
            \Phi_2(c,f_c(z))=\Phi_2(c,z)(z^2-z+c+1).
        \end{equation}
        Substituting $z=f_c^{\circ(m-1)}(0)$ and using
        $f_c^{\circ m}(0)=\bigl(f_c^{\circ(m-1)}(0)\bigr)^2+c$ yields
        \begin{equation}\label{eq:period-two-entry-quotient}
            \frac{\Phi_2(c,f_c^{\circ m}(0))}
            {\Phi_2(c,f_c^{\circ(m-1)}(0))}=H_m(c).
        \end{equation}
        We may evaluate this identity of rational functions at any parameter
        for which the denominator is nonzero. After clearing denominators,
        the polynomial identity holds at every parameter.

        At a root $c$ of the monic polynomial $H_m\in\Z[c]$ of degree
        $2^{m-1}$, put $x=f_c^{\circ(m-1)}(0)$. Then
        $f_c^{\circ m}(0)=x-1$ and $c=x-1-x^2$, giving
        \begin{equation}\label{eq:period-two-cycle-values}
            f_c^{\circ(m+1)}(0)=-x,\qquad
            f_c^{\circ(m+2)}(0)=x-1.
        \end{equation}
        This is a two-cycle, since $x=1/2$ would give $c=-3/4$,
        contrary to the integrality of $c$. If the preperiod is less than
        $m$, then $x$ and $x-1$ are the cycle points. Both are roots of
        $\Phi_2(c,z)$, so they sum to $-1$, giving $x=0$ and $c=-1$.

        Conversely, suppose $f_c$ has exact critical orbit type $(m,2)$.
        The distinct points $f_c^{\circ(m-1)}(0)$ and
        $f_c^{\circ(m+1)}(0)$ have the same image under $f_c$ and hence
        are negatives of one another. The two cycle points sum to $-1$,
        so $H_m(c)=0$.
        At $c=-1$, the critical orbit alternates between $0$ and $-1$, giving
        $H_m(-1)=2$ for even $m$ and $H_m(-1)=0$ for odd $m$.
        Differentiating the critical-orbit recurrence gives the derivative values
        \begin{equation}\label{eq:period-two-periodic-derivatives}
            \left.\frac{d}{dc}f_c^{\circ r}(0)\right|_{c=-1}
            =\begin{cases}
                1,&r\text{ odd},\\
                -1,&r\text{ even},
            \end{cases}
            \qquad r\ge1.
        \end{equation}
        For odd $m$, $H_m'(-1)=2$, so $c+1$ occurs with multiplicity one.
        By the factorization above, the M\"obius product $F_{2,m,2}$ in
        \cite[Definition~1.3]{Goksel} is the quotient in
        \eqref{eq:period-two-entry-quotient}. At $m=1$ it gives
        $F_{2,1,2}=c+1$. The definition divides out this factor exactly
        when $2\mid m-1$, giving \eqref{eq:Gm2}.
        In \cite{BEK} take $D=2$, $k=m$, and $a=c$. For $c\ne0$,
        the map $z\mapsto cz$ conjugates $cz^2+1$ to $f_c(z)$ and
        preserves the exact critical orbit type. The identification of
        $G_{m,2}$ with Goksel's $G_{2,m,2}$ gives irreducibility over $\Q$
        for every $m\ge2$ by \cite[Corollary~1.5]{Goksel}.
    \end{proof}

    \begin{lemma}[Nonreal CM parameters of eventual period two]\label{lem:excludeTwo}
        Let $c$ be nonreal and belong to a CM field. If $f_c$ has exact
        critical orbit type $(m,2)$ with $m\ge2$, then $m$ equals $2$
        and $c$ belongs to $\{i,-i\}$.
    \end{lemma}
    \begin{proof}
        By Theorem~\ref{thm:irred}, $G_{m,2}$ is monic and irreducible
        over $\Q$ and vanishes at $c$, so it is the minimal polynomial
        of $c$. A real root would therefore be a real conjugate of $c$,
        which Lemma~\ref{lem:CMconj} rules out. If $m$ is odd, the polynomial
        has odd degree $2^{m-1}-1$ and hence has a real root.
        For even $m\ge4$, the values of $f_c^{\circ(m-1)}(0)$ at $-2$
        and $-1$ are $2$ and $-1$.
        Thus some $a\in(-2,-1)$ satisfies $f_a^{\circ(m-1)}(0)=0$, and
        \begin{equation}\label{eq:period-two-real-sign-change}
            H_m(-2)=1,\qquad H_m(a)=a+1<0.
        \end{equation}
        Thus $G_{m,2}=H_m$ has a real root for even $m\ge4$ as well.
        The remaining case $m=2$ gives $G_{2,2}(c)=c^2+1$, with roots $\pm i$.
    \end{proof}

    \begin{proof}[Proof of Theorem~\ref{thm:main} in degree $d=2$]
        Let $c$ lie in a CM field $F$ and suppose $f_c$ is PCF.
        If $c$ is real, Lemmas~\ref{lem:CMconj} and~\ref{lem:totreal}
        give $c\in\{0,-1,-2\}$. Otherwise let $(m,n)$ be its exact
        critical orbit type, with $m=0$ or $m\ge2$ by Subsection~\ref{sec:prelim}.

        \smallskip\noindent
        \emph{Case 1. The critical point is periodic ($m=0$).}
        Since $c\ne0$, we have $n\ge2$. Set $v=f_c^{\circ(n-1)}(0)$,
        which is nonzero because $n$ is the exact period. It is an algebraic
        integer because it belongs to $\Z[c]$ and $c$ is integral by
        Lemma~\ref{lem:escape}.
        The relation $f_c^{\circ n}(0)=0$ gives $c=-v^2$. Since $c\ne-2$,
        Lemma~\ref{lem:arch} yields
        \begin{equation}\label{eq:precritical-cm-bound}
            |\sigma(v)|^2=|\sigma(c)|<2
            \qquad(\sigma:F\hookrightarrow\C).
        \end{equation}
        Theorem~\ref{thm:CMgapIntro} gives $v\in\mui$, hence $c\in\mui$.
        By Lemma~\ref{lem:rootunity}, the only nonreal choices are $\pm i$.
        For both, the critical point is strictly preperiodic, so neither
        is possible in this case.

        \smallskip\noindent
        \emph{Case 2. $m\ge2$ and $n\ge3$.}
        This is excluded by Lemma~\ref{lem:excludeHigher}.

        \smallskip\noindent
        \emph{Case 3. $m\ge2$ and $n=1$.}
        The quadratic case of Goksel's irreducibility result
        \cite[Corollary~1.1(i)]{GokselOrbit} follows by the Eisenstein
        argument of Buff, Epstein and Koch \cite{BEK} for the polynomial
        \begin{equation}\label{eq:fixed-cycle-parameter-polynomial}
            A_m(c)=\frac{f_c^{\circ m}(0)+f_c^{\circ(m-1)}(0)}{c}.
        \end{equation}
        The distinct iterates at times $m-1$ and $m$ have equal images,
        hence are negatives of each other. Thus $A_m(c)=0$ as $c\ne0$.
        For $s\ge1$, the polynomial $f_c^{\circ s}(0)$ in $c$ is monic
        of degree $2^{s-1}$ and belongs to $c+c^2\Z[c]$. Induction gives
        \begin{equation}\label{eq:critical-iterate-mod-two}
            f_c^{\circ s}(0)\equiv\sum_{j=0}^{s-1}c^{2^j}\pmod2.
        \end{equation}
        Consequently $A_m\in\Z[c]$ is monic of degree $2^{m-1}-1$,
        has constant coefficient exactly $2$, and satisfies
        \begin{equation}\label{eq:fixed-cycle-eisenstein-congruence}
            A_m(c)\equiv c^{2^{m-1}-1}\pmod2.
        \end{equation}
        Hence $A_m$ is Eisenstein at $2$ and is the minimal polynomial
        of $c$. Its odd degree gives a real root, so $c$ has a real
        conjugate, contrary to Lemma~\ref{lem:CMconj}.

        \smallskip\noindent
        \emph{Case 4. $m\ge2$ and $n=2$.}
        Lemma~\ref{lem:excludeTwo} gives $m=2$ and $c=\pm i$.

        Conversely, the critical orbits of the five remaining parameters are
        \begin{equation}\label{eq:quadratic-classified-orbits}
            \begin{array}{ll}
                c=0&0\mapsto0,\\
                c=-1&0\mapsto-1\mapsto0,\\
                c=-2&0\mapsto-2\mapsto2\mapsto2,\\
                c=i&0\mapsto i\mapsto-1+i\mapsto-i\mapsto-1+i,\\
                c=-i&0\mapsto-i\mapsto-1-i\mapsto i\mapsto-1-i.
            \end{array}
        \end{equation}
        These orbits are finite, and all five parameters lie in the
        CM field $\Q(i)$.
    \end{proof}

    \begin{proof}[Completion of the proofs of Theorem~\ref{thm:main} and Corollary~\ref{cor:ab}]
        Proposition~\ref{prop:high} and the quadratic case prove
        Theorem~\ref{thm:main}. Every listed parameter lies in $\Q^{\ab}$,
        and every element of $\Q^{\ab}$ lies in a CM field by the argument
        following Theorem~\ref{thm:KW}. This proves Corollary~\ref{cor:ab}.
    \end{proof}

    \begin{corollary}[Exact critical orbit types for $f_{2,c}$]\label{cor:types}
        Let $f_{2,c}$ be PCF with $c$ in a CM field. Its exact critical
        orbit type is $(0,1)$ for $c=0$, $(0,2)$ for $c=-1$, $(2,1)$ for
        $c=-2$, and $(2,2)$ for $c=\pm i$. In particular, its eventual
        period is at most two.
    \end{corollary}

    \section*{Declaration of AI use}
    The author initiated and directed the investigation. The author used
    OpenAI's Codex with GPT-6.1 Sol to explore proof strategies, develop
    and check mathematical arguments, and carry out exploratory
    computations. For an earlier draft, the author used a script written
    with Codex to verify finitely many cases in the quadratic classification.
    Codex proposed the use of cyclotomic averages and helped develop the
    argument that replaced those checks. The proofs in this paper do not
    rely on the script.

    The author also used Codex to search the literature, check citations
    and the applicability of cited results, and draft, reorganize, and
    revise the English manuscript. The author used Codex to write the
    Python code that numerically approximates the Multibrot sets and
    plots the classified parameters in Figure~\ref{fig:cm-pcf-multibrot}.
    The author takes full responsibility for the results, proofs,
    references, figure, and final text.

\end{document}